\documentclass[a4paper,12pt,twoside]{article}

\usepackage[left=2cm,right=2cm,top=3cm,bottom=3cm]{geometry}
\usepackage[utf8]{inputenc}
\usepackage[T1]{fontenc} 
\usepackage[english]{babel}
\usepackage{amsmath}
\usepackage{amsthm}
\usepackage{enumitem} 
\usepackage{stmaryrd} 
\usepackage{amssymb,bm}
\usepackage{tikz,tkz-tab} 
\usepackage{multicol} 
\usepackage{parskip}
\usepackage[colorlinks=true, linkcolor=blue, citecolor=red]{hyperref} 
\usepackage{cleveref}
\usepackage{verbatim}
\usepackage{perpage} 
\usepackage[linesnumbered,ruled,vlined,french]{algorithm2e}
\usepackage{fancyhdr} 
\usepackage{titlesec}
\usepackage{setspace}
\usepackage{mathtools}
\usepackage{suffix}

\newcommand{\mycomment}[1]{}

\renewcommand{\mod}[1]{(\mathrm{mod}~#1)}
\WithSuffix\newcommand\psum*[1]{\mathop{\sum\nolimits^{\mathrlap{#1}}}}

\newcommand{\KK}{\mathbb{K}}
\newcommand{\NN}{\mathbb{N}}
\newcommand{\QQ}{\mathbb{Q}}

\newcommand{\F}{\mathcal{F}}

\renewcommand{\O}{\mathcal{O}}

\newcommand\numberthis{\addtocounter{equation}{1}\tag{\thesection.\theequation}}

\newcommand{\textAuthor}{To fill - Author}
\newcommand{\textTitle}{To fill - Title}
\newcommand{\AuthorIs}[1]{\renewcommand{\textAuthor}{#1}}
\newcommand{\TitleIs}[1]{\renewcommand{\textTitle}{#1}}
\newenvironment{entete}[1]
{
  \begin{center}
  {\Large\bf \uppercase{\textTitle}}
  \end{center}
  \vspace{0.5cm}
  \begin{center}
  {\large \uppercase{\textAuthor}}
  \end{center}
  \vspace{0.5cm}
	   {\begin{center}
	     \begin{minipage}{0.8\textwidth}
	     \small{ABSTRACT}.\footnotesize{~#1}\hfill
	     \end{minipage}
	     \end{center}
	      \vspace{0.2cm}
	      }
}

\renewcommand{\theequation}{\arabic{equation}}

\MakePerPage{footnote} 
\counterwithin*{equation}{section}
\labelformat{section}{Section #1}
\labelformat{subsection}{subsection #1}
\titleformat*{\section}{\large\bfseries\centering}
\newcounter{prop}[section]
\newcounter{theo}[section]

\newcounter{lemm}[section]
\newcounter{coro}[section]
\newcounter{fait}[section]
\newcounter{hypo}[section]
\newcounter{conj}[section]
\newcounter{glob}[section]

\newenvironment{proo}[1]{\par\vspace{1em}\noindent  {\it Proof\ifx\relax#1\relax\else~(#1)\fi:}}{\qed\par}

\makeatletter

\newenvironment{theo}[1]{\par\vspace{1em}\noindent
\refstepcounter{glob}
\protected@edef\@currentlabelname{Theorem \arabic{section}.\arabic{glob}}\phantomsection
\noindent \textbf{Theorem \arabic{section}.\arabic{glob}}\ifx\relax#1\relax\else~(#1)\fi:\it}{
\par\vspace{1em}}

\newenvironment{prop}[1]{\par\vspace{1em} \noindent
\refstepcounter{glob}
\protected@edef\@currentlabelname{Proposition \arabic{section}.\arabic{glob}}\phantomsection
\noindent \textbf{Proposition \arabic{section}.\arabic{glob}}\ifx\relax#1\relax\else~(#1)\fi:\it}{\par\vspace{1em}
}

\newenvironment{lemm}[1]{ \par\vspace{1em} \noindent
\refstepcounter{glob}
\protected@edef\@currentlabelname{Lemma \arabic{section}.\arabic{glob}}\phantomsection
\noindent \textbf{Lemma \arabic{section}.\arabic{glob}}\ifx\relax#1\relax\else~(#1)\fi:\it}{\par\vspace{1em}
}

\newenvironment{rema}{\par\vspace{1em} \noindent
\refstepcounter{glob}
\protected@edef\@currentlabelname{Fact \arabic{section}.\arabic{glob}}\phantomsection\noindent\underline{\bf Remark \arabic{section}.\arabic{glob}}:}
	{\par\vspace{1em}
	}
\makeatother

\TitleIs{On abelian cubic fields with large class number} 
\AuthorIs{Jérémy Dousselin} 

\begin{document}
\entete{We investigate the large values of class numbers of cubic fields, showing that one can find arbitrary long sequences of "close" abelian cubic number fields with class numbers as large as possible. We also give a first step toward an explicit lower bound for extreme values of class numbers of abelian cubic fields.}
\section{Introduction}
While working on the theory of cyclotomic fields, Ernst Kummer understood why attempts to prove Fermat's Last Theorem by factorization methods using roots of unity kept failing. In general, the ring generated by those roots does not satisfy the fundamental theorem of arithmetic, in the sense that the unicity of prime factorization is not guaranteed. Later, in 1876, Dedekind introduced the sets named {\it ideals}, in the third edition of {\it Vorlesungen über Zahlentheorie}. The concept of the ideal class group of a ring $R$ is then formalized, and its cardinal $h_R$ -called the class number- is used as a measure of how often the unicity of factorization fails in $R$.

But long before Dedekind's formalism of ideals, it was Gauss who studied what would today be called ideal class groups, in the context of the theory of binary integral quadratic forms. A famous conjecture bears his name, and states that if $h(d)$ is the class number of $\QQ(\sqrt d)$, with $d$ a fundamental discriminant, then $h(d)=1$ for infinitely many positive $d$'s. We know that the case of imaginary quadratic fields is easier to deal with, as shown by Heilbronn \cite{heilbronn} when he proved that $h(d)\to\infty$ when $d\to-\infty$ through the set of negative discriminants. Gauss' conjecture embodies a huge difficulty one faces while studying class numbers of real fields: the presence of non-trivial units heavily affects the size of $h(d)$, making it extremely difficult to understand the behaviour of small $h(d)$'s. The impact of these non-trivial units is captured by the Class Number Formula, which states that 
\[\quad h(d)=\frac{L(1,\chi_d)\sqrt d}{\log \varepsilon_d},\]
for all positive fundamental discriminants $d$, where $\chi_d$ is Kronecker's symbol associated to $d$, and $L(1,\chi_d)$ is the Dirichlet $L$-function attached to $\chi_d$. Finally $\varepsilon_d$ is the fundamental unit, defined by $(t+\sqrt du)/2$ if $(t,u)$ is the smallest solution in $t$ of $t^2-du^2=4$.

In 1977, Montgomery and Weinberger \cite{Montgomery1977} proved that there exists infinitely many real quadratic fields $\QQ(\sqrt d)$ such that 
\[\numberthis\label{montvau}h(d)\gg \sqrt d (\log\log d/\log d).\]
This is widely believed to be the best bound, up to the constant, for it is known that on GRH all fundamental positive discriminants $d$ satisfy
\[\label{h_grh} \numberthis \quad h(d)\leq (4e^\gamma+o(1))\sqrt d\frac{\log\log d}{\log d}.\]
Indeed, assuming GRH Littlewood \cite{littlewood} showed that for any fundamental discriminant $d>0$, we have that 
\[\numberthis\label{littleconj}|L(1,\chi_d)|\leq (2e^\gamma +o(1))\log\log d,\]
where $\gamma$ is the usual Euler-Mascheroni constant. Because $\varepsilon_d>\sqrt d/2$, and hence $\log \varepsilon_d\geq (1/2+o(1))\log d$, the Class Number Formula combined with GRH gives \eqref{h_grh}. More recently, in 2015, Lamzouri \cite{lamzouri} proved there are at least $x^{1/2-1/\log\log x}$ real quadratic fields $\QQ(\sqrt d)$ with discriminant $d\leq x$ such that 
\[\label{ineg-lamzouri}\numberthis h(d)\geq (2e^\gamma+o(1))\sqrt d\frac{\log\log d}{\log d},\]
and that this holds for at most $x^{1/2+o(1)}$ such fields. The constant $e^{2\gamma}$ is widely believed to be best possible, since Littlewood conjectured that one should replace $2e^\gamma$ by $e^\gamma$ in \eqref{littleconj}. Granville and Soundararajan's paper \cite{gvillsound} provides strong support to this conjecture.

\subsection{Large values of class numbers of abelian cubic fields}
In 2004, Duke \cite{duke} investigated a generalization of Montgomery and Weinberger's result \eqref{montvau} to higher degree number fields, and in particular to abelian cubic fields. Note that, in the context of these fields, the Class Number Formula states that for an abelian cubic field $\KK$ of discriminant $d$ and regulator $R$, we have
\[\numberthis\label{eq_class}h_\KK=\frac{d^{1/2}|L(1,\chi)|^2}{4R},\]
for some primitive cubic character $\chi$. Duke's result concerning these fields states that there is an absolute constant $c>0$ such that there are infinitely many abelian cubic fields with arbitrarily large discriminant $d$ for which
\[\label{ine-duke}\numberthis h>cd^{1/2}\left(\frac{\log\log d}{\log d}\right)^2.\]
This bound is easily proved to be best possible, assuming GRH, up to the constant. In his paper, Duke uses specific abelian cubic fields, known as the "simplest cubic fields," which were extensively studied in Shanks' paper \cite{shanks}. These fields are obtained by adjoining to $\QQ$ any root of the polynomial $f_t(x):=x^3-tx^2-(t+3)x-1$, $t\in\NN$ (all three roots of $f_t$ generate the same field, written $\KK_t$). This polynomial is easily shown to have a discriminant given by $\mathrm{disc}(f_t)=g(t)^2$, where $g(t)=t^2+3t+9$. We denote by $\F(x)$ the set of these simplest cubic fields with discriminant $\leq x$.

In the same paper, Duke investigates the case of number fields of degree $n\geq 4$. Assuming both Artin's conjecture and GRH, he is able to prove that there is a constant $c_n>0$ such that there exist totally real number fields $K$ of degree $n$, whose normal closure has the full symmetric group $ S_n$ as its Galois group, with arbitrarily large discriminant $d$ for which
\[h_K>c \sqrt d(\log\log d/\log d)^{n-1}.\]
Recently, in 2020, Lemke Oliver, Thorner and Zaman \cite{Oliver} were able to prove a stronger result unconditionally. They proved that for any fixed integers $r_1,r_2\geq 0$ with $n:=r_1+2r_2\geq 2$, there are number fields $F$ of signature $(r_1,r_2)$ with arbitrarily large discriminant $d$ whose normal closure has $S_n$ as its Galois group, for which
\[h_F\gg_{r_1,r_2} d^{1/2}\left(\frac{\log\log d}{\log d}\right)^{r_1+2r_2-1}.\]

Our first goal in this paper is to extend \eqref{ineg-lamzouri} to abelian cubic fields: more precisely, we want to obtain an explicit constant in \eqref{ine-duke}. In this direction, we prove that
\begin{theo}{}\label{th:2}
Let $x$ be a large real number. For at least $x^{1/4-o(1)}$ fields $\KK$ of $\F(x)$, we have that 
\[h_\KK\geq \left(\frac{4}{91}e^{2\gamma}+o(1)\right)\sqrt d\left(\frac{\log\log d}{\log d}\right)^2,\]
where $d$ is the discriminant of $\KK$.
\end{theo}

This result is to be compared with the upper bound for $h_\KK$ obtained under the assumption of GRH. We would have, for every field $\KK$ of discriminant $d$ in \nameref{th:2}:
\[h_\KK\leq \left(\frac{64}{91}e^{2\gamma}+o(1)\right)\sqrt d\left(\frac{\log\log d}{\log d}\right)^2.\]
This comes from the facts that for all characters $\chi$ of conductor $q$ attached to such a field $\KK$, we have $\chi(2),\chi(3)=e^{\pm 2i\pi/3}$, $\mathrm{reg}(K)=(1/16+o(1))\log^2 d$ (see the proof of \nameref{th:2}), and under GRH we have (see Lemma 2.1 of \cite{gvillsound}):
\[L(1,\chi)=\prod_{p\leq (\log q)^2}\left(1-\frac{\chi(p)}p\right)^{-1}(1+o(1)).\]
Again, it is widely believed that $L(1,\chi)$ can be approximated by the shorter Euler product over the primes $p\leq (\log q)^{1+o(1)}$, and Granville and Soundararajan \cite{gvillsound} gave strong evidence in this direction. Thus, if $\KK$ is a field of discriminant $d$ in \nameref{th:2}, one should rather expect the upper bound
\[h_\KK\leq \left(\frac{16}{91}e^{2\gamma}+o(1)\right)\sqrt d\left(\frac{\log\log d}{\log d}\right)^2.\]

\begin{rema}
Note that the constant in \nameref{th:2} is not the best one can expect, since the cubic case is harder to handle than the quadratic case. The quality of our constant relies on our ability to sieve with precision over certain families of primitive cubic characters. In the context of quadratic characters, Lamzouri \cite{lamzouri} was able to efficiently use  Heath-Brown's quadratic large sieve \cite{heath} to prove \eqref{ineg-lamzouri}. However it seems that cubic large sieve inequalities share an important restriction compared to the quadratic large sieve: they need to be applied to square-free supported numbers. According to Heath-Brown \cite{heath}, this restriction seems like an inevitable technical difficulty if one wants to preserve the sharpness of these sieves. In particular, the Baier and Young cubic large sieve (see Theorem 1.4 of \cite{young}) does not avoid this restriction. While we tried to use this sieve in the method of Granville and Soundararajan (\cite{gvillsound}, proof of Proposition 2.2), the loss caused by the sum being taken only over square-free numbers prevented us from getting a better constant than the one implied by the general large sieve.

Moreover, the family of simplest cubic fields studied by Duke is more complex than Chowla's family of real quadratic fields, which was used by Lamzouri \cite{lamzouri} to study large values of quadratic class numbers. Both these difficulties harmed the quality of our constant in \nameref{th:2}.
\end{rema}
\begin{rema}
Another idea to attack the problem would be to use the fact that the families of characters we have to sieve over have polynomial moduli. The large sieve for characters to polynomial moduli is a recurrent topic in number theory, and many authors such as Baier, Zhao \cite{baierzhao},\cite{zhao}, Munsch \cite{munsch}, and Halupczok \cite{halupczok}, studied it. Unfortunately,  it seems that currently none of these large sieve inequalities to polynomial moduli is sharp enough to yield a better constant in \nameref{th:2} than the one given by the general large sieve. However, a rather strong conjecture about the large sieve to polynomial moduli (see (2.7) of \cite{Mhalu}) may help. If one assumes this conjecture, one could obtain the constant $16/273$ in \nameref{th:2}. This would save a factor of $4/3$.
\end{rema}

\subsection{Tuples of abelian cubic fields with large class numbers, and whose discriminants are close}
Very recently Cherubini, Fazzari, Granville, Kala and Yatsyna \cite{granville} studied the number of consecutive quadratic fields with large values of class numbers. Surprisingly, they were able to prove that for any fixed $k\geq 1$, there are at least $x^{1/2-o(1)}$ integers $d\leq x$ such that 
\[\forall j=1,...,k,\quad h_{\QQ(\sqrt{d+j})}\gg_k\frac{\sqrt d}{\log d}\log\log d.\]

The second goal of this article is to adapt this result to the cubic fields studied by Duke. To extend this to the cubic case, there are multiple ways of translating this "number fields proximity." One can try to prove an analogue for pure cubic fields, considering a sequence 
\[\QQ(\sqrt[3]{d+j}),\quad j=1,...,k.\]
This probably is too hard to prove at the moment, so we instead choose to translate the closeness of the fields by the proximity of their discriminant. More explicitly, we will prove the following:
\begin{theo}{}\label{th:1}
We fix an integer $k\geq 1$ and let $x$ be a large real number. There are $\gg_k x^{1/4-o(1)}$ $k$-tuples of distinct abelian cubic fields, say $(\KK^1,...,\KK^k)$, of discriminants $x\leq D_1\leq...\leq D_k\leq 2x$, such that 
\[\forall j=1,...,k-1:\quad  D_{j+1}-D_j\leq D_j^{3/4+o(1)}\]
and
\[\forall j=1,...,k: \quad h_{\KK^j}\gg_k\sqrt {D_j} \left(\frac{\log \log D_j}{\log D_j}\right)^2,\]
both as $x\to\infty$.
\end{theo}

\noindent {\bf\underline{Acknowledgement:}} I would like to thank my PhD advisor, Youness Lamzouri, for his guidance through the preparation of this article.

\section{Tuples of fields with large class numbers: proof of \nameref{th:1}}
We fix $k\geq 1$ and we let $x$ be a large real number. Thanks to the Class Number Formula \eqref{eq_class}, we know that proving our result relies on our ability to find many simplest cubic fields with large $|L(1,\chi)|$ and small regulator. If we keep all the simplest cubic fields in $\F(x)$, then we might not be able to easily estimate the regulator/discriminant of some of them. Fortunately, the following lemma will highlights one kind of simplest cubic fields that is easy to handle.
\begin{lemm}{Lemma 1 of \cite{duke}}\label{lemma:duke}
Let $t\in\NN$. If $g(t)$ is squarefree, then $D_{t}:=\mathrm{disc}(\KK_{t})=g(t)^2$ and 
\[\mathrm{reg}(\KK_{t})=\frac1{16}(1+o(1))\log^2(D_{t}).\]
\end{lemm}
Note that for any field $\KK_t$ such that $g(t)$ is squarefree, the character $\chi$ appearing in the Class Number Formula is a primitive Dirichlet character of order $3$ and conductor $g(t)$. 

The strategy is now to produce $k$-tuples of the form $(\KK_{t+\delta_1},...,\KK_{t+\delta_k})$, where $\delta_1,...,\delta_k\ll x^{o(1)}$, so that enough squarefree values of $g(t+\delta_1),...,g(t+\delta_k)$ exist when $t$ is restricted to a well chosen arithmetic progression. This arithmetic progression should be such that a lot of tuples have fields with large $|L(1,\chi)|$, but we will deal with this later. 
\begin{lemm}{}\label{lemma:sieve}
Let $x$ be a large real number, let $\varepsilon>0$ be a small real number, let $\delta_1,...,\delta_k$ be integers $\ll x^{2\varepsilon}$, $q=x^{\varepsilon(1+o(1))}$, and let $a$ be an integer such that $t\equiv a\mod q$ implies that $g(t+\delta_1),...,g(t+\delta_k)$ have no prime divisor $\leq \varepsilon\log x$. We fix $\alpha\in[0.02,1]$ a constant. We define $N_\alpha(x;a,q)$ to be the number of $x\leq t\leq(1+ \alpha) x$, $t\equiv a\mod q$ such that $g(t+\delta_1),g(t+\delta_2),...,g(t+\delta_k)$ are all squarefree. Then
\[N_\alpha(x;a,q)\gg_{k} x^{1-2\varepsilon}.\]
\end{lemm}
\begin{proo}{}
By our hypothesis if $t\equiv a\mod q$, then $g(t+\delta_j)$, $j=1,...,k$, cannot have a prime divisor lower than $\varepsilon\log x$.

Now, we want to study the divisibility of $g(t+\delta_j)$, $j=1,...,k$, by larger primes when $t$ is restricted to the arithmetic progression $t\equiv a\mod q$. As detailed through \cite{granville} (see (2.10) and (2.11)), a {\it soupçon} of sieve theory will be the only ingredient we need.
\begin{enumerate}[label=\roman*)]
\item First of all, we put $z=q^2(\log x)^{4k}$, and we want to show that $\varepsilon\log x<p\leq z$ does not divide any $g(t+\delta_j)$, $j=1,...,k$, for $x^{1-o(1)}$ integers $x\leq t\leq (1+ \alpha) x$ in the arithmetic progression $t\equiv a\mod q$. We fix such a $p$. We are studying $k$ quadratic polynomials, and for each of them, there are at most two classes $t\mod p$ such that $g(t+\delta_j)\equiv 0\mod p$. Therefore, there are at most $2k$ classes $t\mod p$ such that at least one $j\in\{1,...,k\}$ satisfies $g(t+\delta_j)\equiv 0\mod p$. Thus, the fundamental theorem of sieve theory ensures that the number of unsieved integers is
\begin{align*}
\numberthis\label{sieve1}\gg \frac{x}q\prod_{\varepsilon\log x<p\leq z}\left(1-\frac{2k}p\right)\gg_k \frac{x}{q(\log z)^{2k}}.
\end{align*}

\item Now, we deal with the case of larger primes. We will show that if $z<p\leq 2(1+ \alpha) x/z^{1/2}$, then $p^2$ does not divide any $g(t+\delta_j)$, $j=1,...,k$, for $x^{1-o(1)}$ integers $x\leq t\leq (1+ \alpha) x$ in the arithmetic progression $t\equiv a\mod q$. We fix such a prime $p$. Note that, for an integer $X$, $g(X)$ is divisible by $p^2$ if and only if $4g(X)=(2X+3)^2+27$ is, which is possible for at most two classes $X\mod{p^2}$. Therefore, there are at most $2k$ classes $t\mod{p^2}$ such that $p^2$ divides at least one of the $g(t+\delta_j)$, and hence the number of such $x\leq t\leq (1+ \alpha)x$ in our arithmetic progression is bounded by
\[2k\left(1+\frac{\alpha x}{qp^2}\right),\]
for each of our $p$'s. Then, the number of $t$ removed this way is 
\[\numberthis\label{sieve2}\ll_k\sum_{z<p\leq 2x/z^{1/2}}2k\left(1+\frac{\alpha x}{qp^2}\right)\ll_k\frac{x}{z^{1/2}\log x}+\frac{x}{qz\log z}.\]
\end{enumerate}
By our choice of $z$, we know that the expression in \eqref{sieve1} is larger than the one in \eqref{sieve2}. Since $q\ll x^{3\varepsilon/2}$, there are $\gg_k x^{1-2\varepsilon}$ integers $x\leq t\leq (1+ \alpha)x$, $t\equiv a\mod q$ such that both previous assertions hold. For any such $t$, we know that if $p^2|g(t+\delta_j)$, for some $j$, then $p>2(1+ \alpha)x/z^{1/2}$. We may write $g(t+\delta_j)=\ell p^2$, and we have
\[\ell=\frac{g(t+\delta_j)}{p^2}\leq \frac{2(1+ \alpha)^2x^2}{(2(1+ \alpha)x/z^{1/2})^2}<z.\]
Since $g(t+\delta_j)$ is supposed to have no prime divisor lower than $z$, we may deduce that $\ell=1$. Therefore, $g(t+\delta_j)=p^2$. Writing $g(t+\delta_j)=((2(t+\delta_j)+3)^2+27)/4$, we are lead to write that 
\[(2(t+\delta_j)+3)^2-4p^2=-27,\]
which is not possible since $p$ is large enough. Indeed, if it were possible, we would be able to factorize $-27$ as the product of two integers with a difference of $4p$, which is absurd. Thus, we have proved the lemma.
\end{proo}

Now, it remains to choose suitable integers $a,q,\delta_1,...,\delta_k$, so that our arithmetic progression $t\equiv a\mod q$ produces tuples of fields with large $|L(1,\chi)|$.
To force $|L(1,\chi)|$ to be large over cubic characters $\chi$ attached to the family of simplest cubic fields, except for a negligible set of characters, one may force $\chi(p)$ to be equal to $1$ for many small primes. The reason behind this, for example, is the following approximation formula for $L(1,\chi)$ given by Granville and Soundararajan \cite{gvillsound}, which they proved using zero-density estimates and the large sieve:
\begin{prop}{Proposition 2.2 of \cite{gvillsound}}\label{formule:L}
Let $A\geq 1$ be fixed. Then, for all but at most $Q^{2/A+o(1)}$ primitive characters $\chi\mod q$ with $q\leq Q$ we have
\[L(1,\chi)=\prod_{p\leq (\log Q)^A}\left(1-\frac{\chi(p)}p\right)^{-1}\left(1+\O\left(\frac{1}{\log\log Q}\right)\right).\]
\end{prop}
Moreover, note that one may force $\chi(p)=1$ by making $p$ completely splits in the cubic field associated with $\chi$. Now, we have everything we need to prove our theorem.
\begin{proo}{of \nameref{th:1}}
We fix a prime $p\geq 3k+2$ and $\varepsilon>0$ small enough. Note that $f_t$ splits into 3 different linear factors for a lot of incongruent $t\mod p$, as recalled through the proof of Theorem 1 in \cite{duke}: precisely, for $(p-4)/3$ such $t$'s if $p\equiv 1\mod3$, and for $(p-2)/3$ such $t$'s if $p\equiv 2\mod3$. Thus, with our choice of $p$, there are at least $k$ distinct $t\mod p$ such that $p$ completely splits in $\KK_t$. Denote them by $t_{p,1},...,t_{p,k}$, and remark that for any $j\in\{1,...,k\}$, we have $g(t_{p,j})\not\equiv 0\mod p$. We define $a_j$, $j=1,...,k$, thanks to the Chinese Remainder Theorem, to be the smallest positive integer such that
\[\left\{\begin{array}{ll}
a_j\equiv 2 \mod{13}&\mathrm{if~}13<3k+2,
\\a_j\equiv 1 \mod p&\mathrm{if~}p< 3k+2\text{ and }p\neq 13,
\\a_j\equiv t_{p,j}\mod p&\mathrm{if~}3k+2\leq p\leq \varepsilon\log x.
\end{array}\right.\]
We define $q:=\prod_{p\leq \varepsilon\log x}p$. With this definition, it is clear that for every $j$, $a_j\leq q= x^{\varepsilon(1+o(1))}$, by the Prime Number Theorem. We also introduce $\delta_j:=a_j-a_1$, and we deduce that
\[\numberthis\label{eq:tnp}\delta_j\ll x^{2\varepsilon}.\]
Also note that $\delta_1=0$. Thus, every prime $3k+2\leq p\leq \varepsilon \log x$ splits completely in any field of the form $\KK_{t+\delta_j}$, $t\equiv a_1\mod q$, $j\in\{1,...,k\}$. Moreover, it is easy to check that for any prime $p\leq \varepsilon\log x$ and any $j=1,...,k$, we have that 
\[\forall t\equiv a_1\mod q,\quad g(t+\delta_j)\equiv g(a_j)\not\equiv 0\mod p.\]

If we choose $a:=a_1$ and fix $\alpha\in[0.02,1]$ to be chosen later, then our \nameref{lemma:sieve} yields the following: there are $\gg_k x^{1/4-\varepsilon/2}$ integers $x^{1/4}\leq t\leq (1+\alpha)x^{1/4}$ in the arithmetic progression $t\equiv a_1\mod q$ such that $g(t+\delta_1)$, $g(t+\delta_2)$, ..., $g(t+\delta_k)$ are all squarefree. Therefore \nameref{lemma:duke} implies that there are $\gg_k x^{1/4-\varepsilon/2}$ $k$-tuples of abelian cubic fields $(\KK_{t+\delta_1},...,\KK_{t+\delta_k})$, each with discriminant $D_j(t)=D_j:=g(t+\delta_j)^2$, with $x^{1/4}\leq t\leq (1+ \alpha) x^{1/4}$, $t\equiv a_1\mod q$, and whose regulators are equal to $\frac{1}{16}(1+o(1))\log^2(D_j)$. Furthermore by construction, every prime $3k+2\leq p\leq \varepsilon\log x$ splits completely in each of these fields, and hence their characters $\chi$ satisfy $\chi(p)=1$, for all $3k+2\leq p\leq \varepsilon\log x$. Taking $A$ large enough in \nameref{formule:L} and using the Class Number Formula \eqref{eq_class}, one concludes that $\gg_k x^{1/4-\varepsilon/2}$ of these fields have class numbers with the expected extreme values.

Using \eqref{eq:tnp}, for all $j$ and all $x^{1/4}\leq t\leq (1+ \alpha) x^{1/4}$, $t\equiv a_1\mod q$, we get that $D_j=g(t+\delta_j)^2\sim t^4$. Therefore, for all $j$ and all such $t$, we have that
\[\frac34 x\leq D_j\leq \frac43(1+ \alpha)^4x.\]
Choosing $\alpha=\sqrt[4]{9/8}-1\approx 0.03$ so that the constant on the right-hand side is $3/2$ and changing our $x$ to $X:=(3/4)x$, we have that for all $j$, $X\leq D_j\leq 2X$. 

Moreover for $x^{1/4}\leq t\leq (1+ \alpha) x^{1/4}$, we have $x^{2\varepsilon}\leq t^{8\varepsilon}$. Therefore, using \eqref{eq:tnp} and the fact that $g(t)=t^2+\O(t)$, we find that 
\[
D_j-D_1=g(t+\delta_j)^2-g(t)^2=\O\left((g(t)^2)^{3/4+2\varepsilon}\right)=\O\left(D_1^{3/4+2\varepsilon}\right),
\]
which immediately implies that for all $j=1,...,k-1$, we have $D_{j+1}-D_j\leq D_j^{3/4+o(1)}$.
\end{proo}

\section{An explicit lower bound: proof of \nameref{th:2}}
We let $Q$ be a large real positive number. As implied by \nameref{lemma:duke} together with the Class Number Formula, it only remains to estimate $L(1,\chi)$ to prove \nameref{th:2}. We want an estimate for $L(1,\chi)$ over a family of characters $\chi$ with conductor $\leq Q$ whose size is $\asymp Q^{1/2}$ (the ones associated with our simplest cubic fields with discriminant $\leq Q^2$). Therefore, discarding a negligible set of discriminants, we may directly use the approximation formula for $L(1,\chi)$ in \nameref{formule:L} with $A=4+\delta$, where $\delta$ is a small positive real number.
\begin{proo}{of \nameref{th:2}} 
Let $\varepsilon>0$ be a small real number, $A=4/(1-16\varepsilon)$, and $x$ be a large real number. Applying \nameref{formule:L} with $Q=\sqrt{2 x}$ ensures that all but at most $x^{1/4-3\varepsilon}$ cubic primitive characters $\chi$ of conductor $q\leq \sqrt{2x}$ are such that
\[\numberthis\label{final-mert}L(1,\chi)=\prod_{p\leq (\log x)^A}\left(1-\frac{\chi(p)}p\right)^{-1}\left(1+\O\left(\frac{1}{\log\log x}\right)\right).\]

By the work done in the previous section, we also know that for at least $x^{1/4-\varepsilon}$ of fields $\KK\in\F(x)$ of discriminant $d$ and of character $\chi$, we have that 
\begin{enumerate}[label=\roman*)]
\item any prime $5\leq p\leq \varepsilon\log x$ completely splits, and hence $\chi(p)=1$;
\item the regulator of $\KK$ is given by $1/16(1+o(1))\log(d)^2$.
\end{enumerate}
Furthermore, for any $t\in\NN$, $g(t)\equiv 1\mod 2$. If $t$ is non-divisible by $3$, then we also have $g(t)\equiv 1\mod 3$. The simplest cubic fields considered here are such that their conductor, of the form $g(t)$, is squarefree, which implies that $t\not\equiv 0\mod 3$. Therefore, for any of the $x^{1/4-\varepsilon}$ fields considered above, we have that $\chi(2),\chi(3)\neq 0$. Similarly, one shows that $\chi(2),\chi(3)\neq 1$, since for every $t\in\NN$, $f_t$ does not split into $3$ distinct linear factors $\mod 2$ nor $\mod 3$. Thus, if $\omega:=e^{2i\pi/3}$, then our fields are such that their associated characters satisfy $\chi(2),\chi(3)\in\{\omega,\overline \omega\}$.

Then, for at least $x^{1/4-2\varepsilon}$ fields $\KK\in\F(x)$ with character $\chi$, i), ii) and \eqref{final-mert} are all true, and Mertens' theorem implies that
\begin{align*}
|L(1,\chi)|&=\left|1-\frac{\chi(2)}2\right|^{-1}\left|1-\frac{\chi(3)}3\right|^{-1}\prod_{5\leq p\leq \varepsilon \log x}\left(1-\frac{1}p\right)^{-1}\prod_{\varepsilon\log x<p\leq (\log x)^A}\left|1-\frac{\chi(p)}p\right|^{-1}(1+o(1))
\\&\geq \frac2{\sqrt{91}} e^\gamma\log\log x\prod_{\varepsilon\log x<p\leq (\log x)^A}\left|1-\frac{\chi(p)}p\right|^{-1}(1+o(1)).
\end{align*}
Again by Merten's theorem, we know that 
\begin{align*}
\prod_{\varepsilon\log x<p\leq (\log x)^A}\left|1-\frac{\chi(p)}p\right|^{-1}&=\exp\left(\Re\sum_{\varepsilon\log x<p\leq (\log x)^A}\frac{\chi(p)}p\right)(1+o(1))
\\&\geq\exp\left(-\frac12\sum_{\varepsilon\log x<p\leq (\log x)^A}\frac{1}{p}\right)(1+o(1))
\\&= \frac{1}{\sqrt A}(1+o(1)).
\end{align*}

Putting these together, we get that for at least $x^{1/4-2\varepsilon}$ fields $\KK\in\F(x)$, of discriminant $d$ and character $\chi$, we have that
\[|L(1,\chi)|\geq \left(\frac{1}{\sqrt{91}}+o(1)\right)e^\gamma\log\log d\]
and $\textrm{reg}(\KK)=(1/16+o(1))(\log d)^2$.

Inserting these estimates in \eqref{eq_class}, one concludes.
\end{proo}

\bibliographystyle{abbrv}
\bibliography{biblio}
\end{document}